\documentclass[reqno]{amsart}
\usepackage[hyphens]{url}
\usepackage{hyperref}
\hypersetup{breaklinks=true}
\usepackage{breakurl}
\usepackage{amssymb}
\usepackage{amsthm}
\usepackage{mathtools}
\usepackage{enumitem}
\usepackage{array}
\usepackage{algorithm}

\newcommand{\ZZ}{\ensuremath{\mathbb{Z}}}
\newcommand{\QQ}{\ensuremath{\mathbb{Q}}}
\newcommand{\RR}{\ensuremath{\mathbb{R}}}

\newcommand{\be}{\begin{equation}}
\newcommand{\ee}{\end{equation}}

\newcommand{\atan}{\operatorname{atan}}
\newcommand{\atanh}{\operatorname{atanh}}
\newcommand{\M}{\mathsf{M}}

\begin{document}

\title{Computing elementary functions using multi-prime argument reduction}

\author{Fredrik Johansson}
\address{Inria Bordeaux, 33405 Talence, France}
\email{fredrik.johansson@gmail.com}


\begin{abstract}
We describe an algorithm for
arbitrary-precision computation of the elementary functions
(exp, log, sin, atan, etc.)
which,
after a cheap precomputation,
gives roughly a factor-two speedup over previous
state-of-the-art algorithms at precision from a few thousand bits up to millions of bits.
Following an idea of Sch\"{o}nhage, we perform
argument reduction using Diophantine combinations of
logarithms of primes; our contribution
is to use a large set of primes instead of a single pair,
aided by a fast algorithm to solve the associated integer relation problem.
We also list new, optimized Machin-like formulas for the necessary
logarithm and arctangent precomputations.
\end{abstract}

\maketitle


\section{Introduction}

There are two
families of competitive algorithms for arbitrary-precision computation
of elementary functions:
the first uses Taylor series together with argument reduction
and needs $O(\M(B) \log^2 B)$ time for $B$-bit precision~\cite{brent1976complexity},
while the second is based on
the arithmetic-geometric mean (AGM) iteration for elliptic integrals
and achieves complexity $O(\M(B) \log B)$~\cite{Brent1976}.\footnote{$\M(B)$ is the complexity of $B$-bit multiplication. We can take $\M(B) = O(B \log B)$~\cite{Harvey2021}.}
Due to constant-factor overheads,
optimized implementations of Taylor series
tend to perform better than the AGM
for practical sizes of $B$, possibly even for $B$ in the billions.

The degree of argument reduction is a crucial tuning
parameter in Taylor series methods.
For example, the standard algorithm
for the exponential function\footnote{The logarithmic
and trigonometric functions have analogous algorithms; alternatively, they can be computed
from the exponential via connection formulas and root-finding for inverses.
For a more comprehensive overview of techniques for elementary function evaluation, see Smith~\cite{Smith1989}, Muller~\cite{Muller2016}, Brent and Zimmermann~\cite{mca}, and Arndt~\cite{arndt2010matters}.}
amounts to choosing a reduction parameter $r \ge 0$
and evaluating
\be
\exp(x) = (\exp(t))^{2^r}, \quad t = x/2^r.
\label{eq:expred}
\ee

If $|x| < 1$, this costs~$r$ squarings plus
the summation of $N \le B / r$ terms
of~the series for $\exp(t)$,
or better $N/2$ terms for $s = \sinh(t)$
(with $\exp(t) = s + \smash{\sqrt{s^2 + 1}}$).
In moderate precision (up to around $B = 10^4$) the series
evaluation costs $O(\sqrt{N})$ multiplications;
for quasilinear complexity as $B \to \infty$, we use the
``bit-burst algorithm'': we write
$\exp(t) = \exp(t_1) \cdot \exp(t_2) \cdots$ where $t_j$ extracts~$2^j$ bits in
the binary expansion of $t$ and evaluate each $\exp(t_j)$ series
using binary splitting.

Asymptotically, $r$ should grow at most logarithmically with $B$,
or the $O(r \M(B))$ time spent on squarings will dominate.
In practice, the best $r$ will be of order 10 to 100 (varying with $B$)
and these $r$ squarings may account for a large fraction of the work to evaluate the function.
This prompts the question: can we reduce the argument to size $2^{-r}$ \emph{without} the cost of $r$ squarings?

The only known solution relies on precomputation.
For example, we need only a single multiplication for $r$-bit
reduction if we have
a precomputed table of $\exp(j/2^r)$, $0 \le j < 2^r$,
or $m$ multiplications with an $m$-partite table of $m 2^{r/m}$ entries.
Tables of this kind are useful up
to a few thousand bits~\cite{Johansson2015elementary},
but they are rarely used at higher precision
since they yield diminishing returns
as the space and precomputation time increases linearly
with $B$ and exponentially with $r$.
Most commonly, arbitrary-precision software will only cache
higher-precision values of the constants $\pi$ and $\log(2)$ computed at runtime,
used for an initial reduction to ensure $|x| < 1$.

\subsection*{Sch\"{o}nhage's method}

In 2006, Sch\"{o}nhage \cite{schonhage2006,schonhage2011} proposed
a method to compute elementary functions using ``diophantine combinations of incommensurable logarithms''
which avoids the problem with large tables.
The idea is as follows: given a real number $x$,
we determine integers $c, d$ such that
\be
x \,\approx\, c \log(2) + d \log(3)
\label{eq:xapprox23}
\ee
within some tolerance $2^{-r}$ (it is a standard result in Diophantine
approximation that such $c, d$ exist for any $r$).
We can then use the argument reduction formula
\be
\exp(x) = \exp(t) \, 2^c 3^d, \quad t = x - c \log(2) - d \log(3).
\label{eq:scho23}
\ee
There is an analogous formula for complex $x$
and for trigonometric functions using Gaussian primes.

The advantage of Sch\"{o}nhage's method is that
we only need to precompute or cache the two constants
$\log(2)$ and $\log(3)$ to high precision
while the rational power product $2^c 3^d$ can be computed on the fly
using binary exponentiation.
If $3^c < 2^B$, this step costs $O(\M(B))$.\footnote{In binary arithmetic,
we only need to evaluate $3^c$ since multiplying
by a power of two is free. This optimization is not a vital ingredient of the algorithm, however.}

Sch\"{o}nhage seems to have considered this method
useful only for $B$ in the range from around 50 to 3000 bits (in his words, ``medium precision'').
The problem is that the coefficients $c, d$ in~\eqref{eq:scho23} grow exponentially with
the desired amount of reduction.
Indeed, solutions with $|t| < 2^{-r}$
will generally have $c, d = O(2^{r/2})$.
It is also not obvious how to
compute the coefficients $c$ and $d$ for
a given $x$; we can use a lookup table for small $r$,
but this retains the exponential scaling problem.

\subsection*{Our contribution}

In this work, we describe a version of Sch\"{o}nhage's algorithm
in which we perform reduction using a basis of $n$ primes,
where~$n$ is arbitrary and in practice may be 10 or more.
The coefficients (power-product exponents) will then only
have magnitude around $O(2^{r/n})$,
allowing much greater reduction
than with a single pair of primes.\footnote{Unfortunately, the only published records of Sch\"{o}nhage's algorithm
are two seminar talk abstracts which are light on details.
The abstracts do mention the possibility of combining three primes instead of a single pair ``for an improved design'',
but there is no hint of a practical algorithm working
with arbitrarily large $n$, $r$ and $B$, which will be presented here.}

Section \ref{sect:reduction} presents an algorithm
for quickly finding an approximating linear combination of several logarithms,
which is a prerequisite for making the method practical.
Section~\ref{sect:algorithm} describes the main algorithm
for elementary functions in more detail.
Section \ref{sect:machin} discusses use of
Machin-like formulas for fast precomputation of logarithms
or arctangents, where we tabulate new optimized
multi-evaluation formulas for special sets of values.

Our implementation results presented in section~\ref{sect:implementation}
show that the new version of Sch\"{o}nhage's algorithm scales remarkably well:
we can quickly reduce the argument
to magnitude $2^{-r}$ where we may have $r \ge 100$ at moderately
high precision (a few thousand bits) and perhaps
$r \ge 500$ at millions of bits.
When $n$ is chosen optimally,
the new algorithm runs roughly twice as fast as
the best previous elementary function implementations
(both Taylor and AGM-based)
for bit precisions $B$ from a few thousand up to millions.
The storage requirements ($nB$ bits) and precomputation time
(on par with one or a few extra function evaluations)
are modest enough that
the method is ideal
as a default
algorithm in arbitrary-precision software
over a large range of precisions.

\subsection*{Historical note}

With the exception of Sch\"{o}nhage's work, we are not aware of any
previous investigations into algorithms of this kind for
arbitrary-precision computation of elementary functions
of real and complex arguments.
However, the underlying idea of exploiting differences between
logarithms of prime numbers
in a computational setting
goes back at least to Briggs' 1624~\emph{Arithmetica logarithmica}~\cite{briggs1624,roegel2010reconstruction}.
Briggs used a version of this trick
when extending tables of logarithms of integers.
We revisit this topic in section~\ref{sect:machin}.


\section{Integer relations}

\label{sect:reduction}

We consider the following \emph{inhomogeneous integer relation problem}:
given real numbers $x$ and $\alpha_1, \ldots, \alpha_n$
and a tolerance $2^{-r}$,
find a vector $(c_1, \ldots, c_n) \in \ZZ^n$
with small coefficients such that
\be
x \,\approx\, c_1 \alpha_1 + \ldots c_n \alpha_n
\label{eq:imhrel}
\ee
with error at most $2^{-r}$.
We assume that the equation $c_1 \alpha_1 + \ldots + c_n \alpha_n = 0$
has no solution over the integers.
In the special case where
$P = \{p_1, \ldots p_n\}$ is a set of prime numbers and $\alpha_i = \log(p_i)$,
solving \eqref{eq:imhrel} will find a $P$-smooth rational approximation
\be
\exp(x) \,\approx\, p_1^{c_1} \cdots p_n^{c_n} \in \QQ
\label{eq:imhprod}
\ee
with small numerator and denominator.

Integer relation problems can be solved
using lattice reduction algorithms like LLL~\cite{lenstra1982factoring,Coh1996}.
However, directly solving
\be
c_0 x + c_1 \alpha_1 + \ldots + c_n \alpha_n \approx 0
\ee
will generally introduce a denominator $c_0 \ne 1$,
requiring
a $c_0$-th root extraction on the right-hand side of \eqref{eq:imhprod}.
In any case, running LLL each time we want to evaluate
an elementary function will be too slow.

Algorithm~\ref{alg:linred} solves
these issues by precomputing solutions to the homogeneous
equation
$c_1 \alpha_1 + \ldots c_n \alpha_n \approx 0$
and using these relations
to solve the inhomogeneous version \eqref{eq:imhrel} through iterated reduction.

\begin{algorithm}

\caption{Approximate $x \in \RR$ to within $2^{-r}$ by a linear combination
$$x \,\approx\, c_1 \alpha_1 + \ldots + c_n \alpha_n, \quad c_i \in \ZZ$$
given $\alpha_i \in \RR$ which are linearly independent over $\QQ$.
Alternatively, find a good approximation subject to some size constraint $f(c_1, \ldots, c_n) \le M$.}

\label{alg:linred}

\begin{enumerate}[leftmargin=0.65cm]

\item Precomputation (independent of $x$): choose a real convergence factor $C > 1$. For $i = 1, 2, \ldots$, LLL-reduce
$$
\renewcommand\arraystretch{1.1}
\begin{pmatrix} 1 & 0 & \ldots & 0 & \lfloor C^i \alpha_1 + \tfrac{1}{2} \rfloor) \\
                  0 & 1 & \ldots & 0 & \lfloor C^i \alpha_2 + \tfrac{1}{2} \rfloor) \\
                  \vdots & \vdots  & \ddots & \vdots &     \vdots                     \\
                  0 & 0 & \ldots & 1 & \lfloor C^i \alpha_n + \tfrac{1}{2}  \rfloor) \end{pmatrix}.$$
This yields an approximate integer relation
\be
\varepsilon_i = d_{i,1} \alpha_1 + \ldots d_{i,n} \alpha_n, \quad \varepsilon_i = O(C^{-i}).
\label{eq:approxrel}
\ee
(In fact, it yields $n$ such relations; we can choose any one of them.)
We store tables of the coefficients $d_{i,j}$ and floating-point approximations
of the errors $\varepsilon_i$.
We stop after the first $i$ where $|\varepsilon_i| < 2^{-r}$.

\vskip5pt

\item Reduction (given $x$).
\begin{itemize}
\item Let $(c_1, \ldots, c_n) = (0, \ldots, 0)$.
\item For $i = 1, 2, \ldots$, compute $m_i = \lfloor x / \varepsilon_i + 1/2 \rfloor$ and update:
$$(c_1, \ldots, c_n) \; \gets \; (c_1 + m_i d_{i,1}, \ldots, c_n + m_i d_{i,n}),$$
$$x \; \gets \; x - m_i \varepsilon_i.$$
Stop and return the relation $(c_1, \ldots, c_n)$
when $|x| < 2^{-r}$ or when
the next update will give $f(c_1, \ldots, c_n) > M$.
\end{itemize}

\end{enumerate}

\end{algorithm}

\subsection*{Analysis of Algorithm~\ref{alg:linred}}

We assume heuristically that each step in the precomputation phase (1)
succeeds to find a relation~\eqref{eq:approxrel}
with $\varepsilon_i$ within a small factor of $\pm C^{-i}$
and with coefficients $(d_{i,1}, \ldots, d_{i,n})$
of magnitude $O(C^{i/n})$.
We will simply observe that this always seems to be the case in practice;
a rigorous justification would require
further analysis.

It can happen that picking the first integer relation
computed by LLL yields the same relation consecutively ($\varepsilon_i = \varepsilon_{i+1}$).
In that case, we can just pick a different relation
(while keeping the $\varepsilon_i$ sorted)
or skip the duplicate relation. However, a decrease by much
more than a factor $C$ between successive step should be avoided
as it will result in larger output coefficients.

Phase (1) terminates when $i = N \approx r \log(2)/\log(C)$.
The multiplier $m_i$ computed in each step of the phase (2) reduction
has magnitude around $C$.
The coefficients $(c_1, \ldots, c_n)$
at the end of phase (2) will therefore have magnitude around
\be
\sum_{i=1}^{N} C^{i/n+1} = \frac{C^{1/n+1}}{C^{1/n}-1} \left(2^{r/n}-1\right) \approx \frac{C n}{\log(C)} 2^{r/n}
\ee
or perhaps a bit smaller than this since on average
there can be some cancellation.

The prefactor $C/\log(C)$ is minimized when $C = e$,
or in other words it is
theoretically optimal to force $\varepsilon_i = \Theta(\exp(-i))$.
However, this prefactor does not vary strongly with $C$,
and a choice like $C = 2$ (one bit per step)
or $C = 10$ (one decimal per step) may be convenient.

Step $i$ of phase (1) requires LLL-reducing a matrix
with $\beta$-bit entries where $\beta = O(i)$.
The standard complexity bound for LLL is $O(n^{5+\varepsilon} \beta^{2+\varepsilon})$,
so phase (2) costs $O(n^{5+\varepsilon} r^{3+\varepsilon})$.\footnote{The factor $r^{3+\varepsilon}$ can be improved
to $r^{2+\varepsilon}$ using a quasilinear version of LLL~\cite{novocin2011lll}.}
In our application, the tables generated in phase (1) are small
(a few kilobytes) and do not need to be generated at runtime,
so it suffices to note that the computations are feasible for
ranges of $n$ and $r$ of interest; for empirical results, see section~\ref{sect:implementation}.

Phase (2) requires $O(n r)$ arithmetic operations
with $r$-bit numbers,
for a running time of $O(n r^{2+\varepsilon})$.
It is convenient to treat $x$ and $\alpha_i$ as fixed-point numbers
with $r$-bit fractional part.
As an optimization, we can work with a machine-precision
(53-bit) floating-point approximations of $x'$ and the errors $\varepsilon_i$.
We periodically
recompute
$$x' = x - (c_1 \alpha_1 + \ldots + c_n \alpha_n)$$
accurately
from the full-precision values
only when this approximation runs out of precision, essentially every $53/\log_2(C)$ steps.
The resulting algorithm has very low overhead.
We will not consider asymptotic complexity improvements
since $r$ will be moderate (a small multiple of the word size) in our application.

\section{Computation of elementary functions}

\label{sect:algorithm}

Given $x \in \RR$ and a set of prime numbers $P = \{p_1, p_2, \ldots, p_n\}$,
the algorithm described in the previous section allows us to
find integers $c_1, \ldots, c_n$ such that
\be
t = x - (c_1 \log(p_1) + \ldots + c_n \log(p_n))
\ee
is small, after which we can evaluate the real exponential function as
\be
\exp(x) = \exp(t) \, p_1^{c_1} \cdots p_n^{c_n}.
\ee

Algorithm~\ref{alg:exp} describes the
procedure in some more detail.

\begin{algorithm}

\caption{Computation of $\exp(x)$ for $x \in \mathbb{R}$ to $B$-bit precision using argument reduction by precomputed logarithms of primes.}

\label{alg:exp}

\begin{enumerate}[leftmargin=0.65cm]

\item Precomputation (independent of $x$): select a set of prime numbers $P = \{p_1, p_2, \ldots, p_n\}$ with $p_1 = 2$. Compute $\log(p_1), \ldots, \log(p_n)$ to $B$-bit precision.

\item Using Algorithm~\ref{alg:linred}, find an integer relation $x \approx c_1 \log(p_1) + \ldots + c_n \log(p_n)$, attempting to make the error as
small as possible subject to $\| c_1, \ldots, c_n \|_P \le B$. This step can use low precision (about $r$ bits where $2^{-r}$ is the target reduction, in practice no more than a few machine words).

\item Compute the power product $v/w = p_2^{c_2} \cdots p_n^{c_n}$ as an exact
fraction, using binary splitting to recursively split the set of primes in half and using binary exponentiation to compute the individual powers.

\item Calculate $t = x - (c_1 \log(p_1) + \ldots + c_n \log(p_n))$ using the precomputed logarithms.

\item Compute $u = \exp(t)$ using Taylor series: depending on $B$, either use rectangular splitting for the sinh series or use the bit-burst decomposition $t = t_1 + t_2 + \ldots$ with binary splitting (see~e.g.\ \cite{mca} for details).

\item Return $2^{c_1} u v / w$.

\end{enumerate}

\end{algorithm}

\subsubsection*{Remarks}

The bottleneck in the argument reduction
is the cost of evaluating the
power product $p_1^{c_1} \cdots p_n^{c_n} \in \QQ$.
How large coefficients (exponents) should we allow?
A reasonable heuristic, implemented in~Algorithm~\ref{alg:exp}, is to choose coefficients such that the weighted norm
\be
\nu = \| c_1, \ldots, c_n \|_P = \sum_{\substack{i=1 \\ p_i \ne 2}}^n |c_i| \log_2(p_i)
\ee
is smaller than $B$: this ensures that
the rational power product $p_1^{c_1} \cdots p_n^{c_n}$
has numerator and denominator bounded by $B$ bits.
We discount the prime 2 in the norm
with the assumption that we factor out
powers of two when performing binary arithmetic.
If $|x| > 1$, we should use $\log(2)$
alone for the first reduction in Algorithm~\ref{alg:linred} so that the corresponding
exponentiation is free.

We note that when computing the power product, there is
no need to compute GCDs since the numerator
and denominator are coprime by construction.

There is not much to say about numerical issues;
essentially, we need about $\log_2 (\sum_i |c_i| \log(p_i))$ guard
bits to compensate for cancellation in the subtraction,
which in practice always will be less than one extra machine word.
If $|x| \gg 1$, we need an additional $\log_2 |x|$ guard bits
for the accurate removal of $\log(2)$.


\subsection{Numerical example}

We illustrate computing $\exp(x)$ to 10000 digits (or $B = 33220$ bits)
where
$x = \sqrt{2} - 1$,
using $n = 13$ primes.

The following Pari/GP output effectively shows
the precomputations of phase~(1) of Algorithm~\ref{alg:linred} with
convergence rate $C = 10$.
Since $2^{-100} \approx 7.9 \cdot 10^{-31}$, reducing by
32 relations with $C = 10$ is equivalent to $r = 100$
squarings in \eqref{eq:expred}.

\begin{small}
\begin{verbatim}
? n=13; for(i=1, 32, localprec(i+10); P=vector(n,k,log(prime(k)));
  d=lindep(P,i)~; printf("%s %.5g\n", d, d * P~))
[0, 0, 0, 0, -1, 1, 0, 0, 0, 0, 0, 0, 0] 0.16705
[0, 0, 1, 0, -1, 0, -1, 0, 0, 0, 0, 1, 0] -0.010753
[-1, 0, 0, 0, 0, -1, 1, -1, 0, 1, 0, 0, 0] -0.0020263
[-1, 0, 0, 0, -1, 0, 1, -1, 1, -1, 1, 0, 0] -8.2498 e-5
[1, 0, 1, -1, 0, 1, -1, 1, -1, 0, 0, -1, 1] 9.8746 e-6
[0, 1, 0, -1, -1, 0, 2, -1, 0, -1, -1, 1, 1] 1.5206 e-6
[1, -1, 0, 1, 1, 2, -1, 0, -2, 1, -1, -1, 1] 3.2315 e-8
[1, -1, 0, 1, 1, 2, -1, 0, -2, 1, -1, -1, 1] 3.2315 e-8
[1, 0, 4, -1, -2, 0, 0, 2, 0, -2, -2, 1, 1] 4.3825 e-9
[0, -2, 0, 0, -2, 0, 0, 2, -4, 4, -1, 1, 0] -2.1170 e-10
[1, 1, 4, 1, -1, 1, -2, -3, 0, -4, 3, 1, 1] -7.0743 e-11
[0, -2, -1, 0, 2, 4, 4, 0, 3, 1, -6, -1, -3] 3.3304 e-12
[3, 2, -1, -6, 2, 3, -2, -2, 3, 1, 5, -4, -2] 2.5427 e-13
[-4, -2, 4, -4, 3, 1, 7, 0, -3, -4, 4, -7, 3] -9.9309 e-14
[1, -1, -7, -2, 5, 5, -6, 2, 0, -10, 5, 2, 3] -9.5171 e-15
[3, -2, -7, -9, 6, 6, 3, 9, 1, 8, -15, -4, 0] 6.8069 e-16
[-1, 13, -5, -7, -3, -3, -13, 3, 0, -1, 6, -3, 12] -7.1895 e-17
[-2, 3, -2, 2, -15, 16, 4, -7, 11, -15, 0, 9, -4] 8.1931 e-18
[2, 0, -9, -11, -5, -11, 21, 9, -9, -4, -1, -4, 13] 5.6466 e-19
[6, -9, 0, 9, 9, -2, -4, -22, 4, -7, 0, 5, 11] 4.6712 e-19
[1, -27, 22, -14, -2, 0, 0, -27, -3, -5, 18, 10, 9] -1.0084 e-20
[1, 41, -2, 5, -42, 6, -2, 13, 5, 3, -5, 7, -9] -1.3284 e-21
[4, -5, 8, -8, 6, -25, -38, -16, 24, 13, -10, 10, 24] -8.5139 e-23
[4, -5, 8, -8, 6, -25, -38, -16, 24, 13, -10, 10, 24] -8.5139 e-23
[-43, -2, 4, 9, 19, -26, 92, -30, -6, -24, 11, -4, -18] -4.8807 e-24
[8, 38, -4, 34, -31, 60, -75, 31, 44, -32, -1, -43, 17] 2.7073 e-25
[48, -31, 21, -27, 34, -23, -29, 41, -50, -65, 33, 20, 40] 5.2061 e-26
[-41, 8, 67, -84, 7, -22, -58, -35, 17, 58, -18, 13, 40] -7.9680 e-27
[20, 15, 50, -1, 48, 72, -67, -96, 75, 48, -38, -126, 68] 2.7161 e-28
[26, 20, -35, 16, -1, 75, -13, 2, -128, -100, 130, 46, -13] -3.3314 e-29
[-26, -20, 35, -16, 1, -75, 13, -2, 128, 100, -130, -46, 13] 3.3314 e-29
[137, -26, 127, 45, -14, -73, -66, -166, 71, 76, 122, -154, 53] -1.4227 e-31
\end{verbatim}
\end{small}

We prepend the relation $[1, 0, \ldots]$
for an initial reduction by $\log(2)$, and we can eliminate the
duplicate entries.

The phase (2) reduction in Algorithm~\ref{alg:linred} with $x = \sqrt{2}-1$
now yields
the relation
$$[-274, -414, -187, -314, -211, 651, -392, 463, -36, -369, -231, 634, 0]$$
or
$$\exp(x) \approx \frac{2^{c_1} v}{w} = \frac{13^{651} \cdot 19^{463} \cdot 37^{634}}{2^{274} \cdot 3^{414} \cdot 5^{187} \cdot 7^{314} \cdot 11^{211} \cdot 17^{392} \cdot 23^{36} \cdot 29^{369} \cdot 31^{231}}$$
where the numerator and denominator have 7679 and 7678 bits,
comfortably smaller than $B$.

We compute the reduced argument $t = x - \log(2^{c_1} v/w) \approx -1.57 \cdot 10^{-32}$
by subtracting a linear combination of precomputed logarithms.
Now taking 148 terms of the Taylor series for $\sinh(t)$ yields an error smaller than $10^{-10000}$.
Evaluating this Taylor series using rectangular splitting costs roughly $2 \sqrt{148} \approx 24$ full 10000-digit multiplications,
and this makes up the bulk of the time in the $\exp(x)$ evaluation.

For comparison, computing $\exp(x)$ using \eqref{eq:expred} without
precomputation, it is optimal to perform $r \approx 20$ squarings
after which we need 555 terms of the sinh series,
for a cost of $r + 2 \sqrt{555} \approx 67$ multiplications.\footnote{This estimate is not completely accurate because a squaring is somewhat cheaper than a multiplication (theoretically requiring 2/3 as much work).
The same remark also concerns series evaluation, where some operations are squarings.
We also mention that computing $\exp(x/2^r)$ with the bit-burst algorithm
might be faster than using the sinh series at this level of precision,
though probably not by much; we use the sinh series here for the purposes of illustration since the analysis is simpler.}
Alternatively, computing $\log(x)$ using the AGM requires 25 iterations,
where each iteration $a_{n+1}, b_{n+1} = (a_n+b_n)/2, \sqrt{a_n b_n}$
costs at least as much as two multiplications.

Counting arithmetic operations alone, we can thus expect
Algorithm~\ref{alg:exp} to be at least twice as fast as either method
in this example.
As we will see in section~\ref{sect:implementation}, this back-of-the-envelope
estimate is quite accurate.

\subsection{Trigonometric functions}

We can compute the real trigonometric functions
via the exponential function of a pure imaginary argument,
using Gaussian primes $a + bi \in \mathbb{Z}[i]$
for reduction.
Enumerated in order of norm $a^2+b^2$,
the nonreal Gaussian primes are
\be
\label{eq:gaussprimes}
1+i, \, 2+i, \, 3+2i, \, 4+i, \, 5+2i, \, 6+i, \, 5+4i, \, 7+2i, \, 6+5i, \ldots
\ee
where we have discarded entries that are equivalent under
conjugation, negation or transposition of real and imaginary parts
(we choose here, arbitrarily, the representatives in the
first quadrant and with $a \ge b$).

The role of the logarithms $\log(p)$
is now assumed by the \emph{irreducible angles}
\be
\alpha = \frac{1}{i} \left[ \log(a + bi) - \log(a - bi) \right] = 2 \atan\!\left(\frac{b}{a}\right)
\ee
which define rotations
by $e^{i\alpha} = (a+bi)/(a-bi)$
on the unit circle.
We have the argument reduction formula
\be
\cos(x) + i \sin(x) = \exp(i x) = \exp(i (x - c \alpha)) \frac{(a+bi)^c}{(a-bi)^c}, \quad c \in \ZZ
\ee
which can be iterated over a combination of Gaussian primes.
Algorithm~\ref{alg:expi}
computes $\cos(x)$ and $\sin(x)$ together
using this method.

\begin{algorithm}

\caption{Computation of $\cos(x) + i \sin(x) = \exp(ix)$ for $x \in \mathbb{R}$ to $B$-bit precision
using argument reduction by precomputed irreducible angles.}

\label{alg:expi}

\begin{enumerate}[leftmargin=0.65cm]

\item Precomputation (independent of $x$): select a set of Gaussian prime numbers $Q = \{a_1 + b_1 i, \ldots, a_n + b_ni\}$ from \eqref{eq:gaussprimes} with $a_1+b_1 = 1+i$.
Compute $2 \atan(b_1/a_1), \ldots, 2 \atan(b_n/a_n)$ to $B$-bit precision.

\item Using Algorithm~\ref{alg:linred}, find an integer relation $x \approx c_1 2 \atan(b_1/a_1) + \ldots + c_n 2 \atan(b_n/a_n)$, attempting to make the error as
small as possible subject to $\| c_1, \ldots, c_n \|_Q \le B$. This step can use low precision (about $r$ bits where $2^{-r}$ is the target reduction, in practice no more than a few machine words).

\item Compute the power product
\be
\frac{v}{w} = \frac{(a_2+b_2 i)^{c_2} \cdots (a_n+b_n i)^{c_n}}{(a_2-b_2 i)^{c_2} \cdots (a_n-b_n i)^{c_n}} \in \mathbb{Q}(i)
\label{eq:gaussprod}
\ee
as an exact
fraction, using binary splitting to recursively split the set of primes in half and using binary exponentiation to compute the individual powers.

\item Calculate $t = x - (c_1 2 \atan(b_1/a_1) + \ldots + c_n 2 \atan(b_1/a_1))$ using the precomputed arctangents.

\item Compute $u = \exp(i t)$ using Taylor series (depending on $B$, either using rectangular splitting for the sin series or using the bit-burst decomposition $t = t_1 + t_2 + \ldots$ with binary splitting).

\item Return $i^{c_1} u v / w$.

\end{enumerate}

\end{algorithm}

\subsubsection*{Remarks}

Here, a suitable norm is
\be
\nu = \| c_1, \ldots, c_n \|_Q = \sum_{\substack{j=1 \\ p_j \ne 1+i}}^n |c_j| \log_2(a_j^2+b_j^2).
\ee

The special prime 2 in the argument reduction for the real exponential
is here replaced by the Gaussian prime $1+i$,
for which
\be
\frac{(1+i)^c}{(1-i)^c} = i^c
\ee
can be evaluated in constant time;
the angle reduction corresponds to removal of multiples of $\pi/2$.

We only need to compute the factors in the numerator of the right-hand
side of \eqref{eq:gaussprod} since the remaining product
can be obtained via complex conjugation.
As in the real case, all factors are coprime so we can multiply
numerators and denominators using arithmetic in $\ZZ[i]$ without the need for GCDs.

We can save a marginal amount of work (essentially
in the last division) if we want
either the sine of the cosine alone,
or if we want $\tan(x)$.

\subsection{Inverse functions}

The formulas above can be transposed to compute the
inverse functions. For example,
\be
\log(x) = \log\left(\frac{x}{p_1^{c_1} \cdots p_n^{c_n}}\right) + (c_1 \log(p_1) + \ldots + c_n \log(p_n)).
\ee

For the complex logarithm or arctangent, we need to be careful
about selecting the correct branches.

As an alternative, we recall the standard method
of implementing the inverse
functions using Newton iteration,
starting from an low-precision approximation obtained with any
other algorithm.
The constant-factor overhead of Newton iteration can be reduced
with an $m$-th order method derived
from the addition formula for the exponential function~\cite[section 32.1]{arndt2010matters}.
If $y = \log(x) + \varepsilon$, then
\be
\log(x) = y + \log(1+\delta), \quad \delta = x \exp(-y) - 1.
\label{eq:logy}
\ee

We first compute $y \approx \log(x)$ at precision $B/m$ (calling the same
algorithm recursively until we hit the basecase range) so
that the unknown error $\varepsilon$ is $O(2^{-B/m})$.
Then, we evaluate \eqref{eq:logy} at precision $B$
using the Taylor series for $\log(1+\delta)$ truncated to order $O(\delta^m)$.
This gives us $\log(x)$ with error $O(2^{-B})$.

The inverse trigonometric functions can
be computed analogously via the arc\-tangent:
if $y = \atan(x) + \varepsilon$, then
\be
\atan(x) = y + \atan(\delta), \quad \delta = \frac{x-t}{1 + t x} = \frac{c x - s}{c + s x}, \quad t = \tan(y) = \frac{s}{c} = \frac{\sin(y)}{\cos(y)}.
\label{eq:atany}
\ee

With a suitably chosen $m$ (between 5 and 15, say)
and rectangular splitting for the short Taylor series evaluation,
the inverse functions are perhaps 10\%-30\%
more expensive than the forward functions with this method.

\section{Precomputation of logarithms and arctangents}

\label{sect:machin}

The precomputation of logarithms and arctangents of small integer
or rational arguments is best
done using binary splitting evaluation
of trigonometric and hyperbolic arctangent series
\be
\operatorname{atan}\!\left(\frac{1}{x}\right) = \sum_{k=0}^{\infty} \frac{(-1)^k}{(2k+1)} \frac{1}{x^{2k+1}}, \quad \atanh\!\left(\frac{1}{x}\right) = \sum_{k=0}^{\infty} \frac{1}{(2k+1)} \frac{1}{x^{2k+1}}.
\label{eq:atanseries}
\ee

We want the arguments $x$ in \eqref{eq:atanseries} to be integers,
and ideally large integers so that the series converge rapidly.
It is not a good idea to use the primes~$p$ or
Gaussian integer tangents $b/a$ directly as input since
convergence will be slow;
it is better to recycle values and evaluate
differences of arguments (Briggs' method).
For example, if we have already computed $\log(2)$,
we can compute logarithms of successive primes
using~\cite{gourdon2004logarithmic}
\be
\log(p) = \log(2) + \frac{1}{2} \left( \log\!\left(\frac{p-1}{2}\right) + \log\!\left(\frac{p+1}{2}\right)\right) + \atanh\!\left(\frac{1}{2p^2-1}\right).
\label{eq:logpdiff}
\ee

Methods to reduce arctangents to sums
of more rapidly convergent
arctangent series
have
been studied by Gauss, Lehmer, Todd and others~\cite{lehmer1938arccotangent,todd1949problem,wetherfield1996enhancement}.
The prototype is Machin's formula
\be
\frac{\pi}{4} = \atan(1) = 4 \atan\!\left(\frac{1}{5}\right) - \atan\!\left(\frac{1}{239}\right).
\label{eq:machin}
\ee

\subsection{Simultaneous Machin-like formulas}

If we have the option of computing the set of values
$\log(p_1), \ldots, \log(p_n)$
or $\atan(b_1/a_1), \ldots, \atan(b_n/a_n)$
in any order (not necessarily one by one),
then we can try to
look for optimized simultaneous Machin-like formulas~\cite{arndt2010matters}.

Given the first $n$ primes, we will thus look
for a set of integers $X = \{x_1, x_2, \ldots, x_n\}$,
as large as possible,
such that there is an integer relation
\be
\begin{pmatrix} \log(p_1) \\ \vdots \\ \log(p_n) \end{pmatrix} = M \begin{pmatrix} 2 \atanh(1/x_1) \\ \vdots \\ 2 \atanh(1/x_n) \end{pmatrix}, \quad M \in \QQ_{n \times n}
\label{eq:intrel1}
\ee
or similarly (with different $X$ and $M$) for Gaussian primes
\be
\begin{pmatrix} \atan(b_1/a_1) \\ \vdots \\ \atan(b_n/a_n) \end{pmatrix} = M \begin{pmatrix} \atan(1/x_1) \\ \vdots \\ \atan(1/x_n) \end{pmatrix}, \quad M \in \QQ_{n \times n}.
\label{eq:intrel2}
\ee

For example,
the primes $P = \{2, 3\}$ admit the simultaneous Machin-like formulas
$\log(2) = 4 \atanh(1/7) + 2 \atanh(1/17)$, $\log(3) = 6 \atanh(1/7) + 4 \atanh(1/17)$, i.e.\
$$X = \{7, 17\}, \quad M = \small \begin{pmatrix} 2 & 1 \\ 3 & 2 \end{pmatrix}.$$

The following method to find relations goes back to Gauss who
used it to search for generalizations of Machin's formula.
Arndt~\cite[section 32.4]{arndt2010matters} also discusses
the application of simultaneous computation
of logarithms of several primes.

The search space for candidate sets $X$ in \eqref{eq:intrel1} and \eqref{eq:intrel2}
is a priori infinite,
but it can be narrowed down as follows.
Let $P = \{p_1, \ldots, p_n\}$.
Since $$2 \atanh(1/x) = \log(x+1) - \log(x-1) = \log\!\left(\frac{x+1}{x-1}\right),$$
we try to write each $p \in P$ as a power-product of
$P$-smooth rational numbers of the form $(x+1)/(x-1)$.
We will thus look for solutions $X$ of \eqref{eq:intrel2} of the form
\be
X \subseteq Y, \quad Y = \{ x: x^2 -1 \text{ is } P\text{-smooth}\},
\label{eq:setY}
\ee
i.e.\ such that both $x+1$ and $x-1$ are $P$-smooth.
Similarly, we look for solutions of \eqref{eq:intrel2} of the form
\be
X \subseteq Z, \quad Z = \{ x: x^2 + 1 \text{ is } Q\text{-smooth}\}
\label{eq:setZ}
\ee
where $Q$ is the set of norms $\{a_1^2+a_1^2, \ldots, a_n^2+b_n^2\}$.

It is a nontrivial fact that the sets $Y$ and $Z$ are finite for each fixed set of primes $P$ or $Q$.
For the 25 first primes $p < 100$, the set $Y$ has 16223 elements
which have been tabulated by Luca and Najman~\cite{Luca2010,Luca2013};
the largest element\footnote{Knowing this upper bound,
the Luca-Najman table can be reproduced with a brute force enumeration
of 97-smooth numbers $x-1$ with $x \le 19182937474703818751$,
during which one saves the values $x$ for which trial division shows that $x+1$ is 97-smooth.
This computation takes two hours on a 2022-era laptop.
Reproducing the table $Z$ takes one minute.} is $x = 19182937474703818751$
with
$$x - 1 = 2 \cdot 5^{5} \cdot 11 \cdot 19 \cdot 23^{2} \cdot 29 \cdot 59^{4} \cdot 79,$$
$$x + 1 = 2^{22} \cdot 3 \cdot 17^{3} \cdot 37 \cdot 41 \cdot 43 \cdot 67 \cdot 71.$$
For the first 22 Gaussian primes, having norms $a^2+b^2 < 100$, the set $Z$ has 811 elements
which have been tabulated by Najman~\cite{najman2010smooth};
the largest element is $x = 69971515635443$
with
$$x^2 + 1 = 2 \cdot 5^{5} \cdot 17 \cdot 37 \cdot 41^{2} \cdot 53^{2} \cdot 89 \cdot 97^{3} \cdot 137^{2} \cdot 173.$$

Given a candidate superset $Y = \{ y_1, \ldots, y_r \}$ or $Z = \{ z_1, \ldots, z_s \}$, we can find
a formula $X$ with large entries using linear algebra:

\begin{itemize}
\item Let $X = \{\}$, and let $R$ be an initially empty ($0 \times n$) matrix.
\item For $x = y_r, y_{r-1}$, $\ldots$ or $x = z_r, z_{s-1}$, $\ldots$ in order of decreasing magnitude, let $E = (e_1, \ldots, e_n)$
be the vector of exponents in the factorization of the
rational number $$(x+1)/(x-1) = p_1^{e_1} \cdots p_n^{e_n},$$
respectively,
$$x^2+1 = (a_1^2+b_1^2)^{e_1}, \ldots, (a_n^2+b_n^2)^{e_n}.$$
\item If $E$ is linearly independent of the rows of $R$,
add $x$ to $X$ and
adjoin the row $E$ to the top of $R$; otherwise continue with the next candidate $x$.
\item When $R$ has $n$ linearly independent rows, we have found a complete basis~$X$ and
the relation matrix is given by $M = R^{-1}$.
\end{itemize}

Tables~\ref{tab:logrelations} and ~\ref{tab:atanrelations} give the
Machin-like formulas found with this method
using the exhaustive Luca-Najman tables for $Y$ and $Z$.
We list only the set~$X$ since the matrix~$M$
is easy to recover with linear
algebra (in fact, we can recover it using LLL without
performing any factorization).
The corresponding
\emph{Lehmer measure} $\mu(X) = \sum_{x \in X} 1 / \log_{10}(|x|)$
gives an estimate of efficiency (lower is better).

\subsection{Remarks about the tables}

We conjecture that the formulas
in Tables~\ref{tab:logrelations} and ~\ref{tab:atanrelations} 
are the best possible (in the Lehmer sense) $n$-term formulas for the respective sets of $n$ primes or Gaussian primes.

Apart from the first few entries which are well known, we are not aware of a
previous tabulation of this kind.
There is an extensive literature about Machin-like formulas
for computing $\pi$ alone, but
little about computing several arctangents simultaneously.
There are some preexisting tables for logarithms,
but they are not optimal. Arndt~\cite{arndt2010matters} gives a slightly less efficient
formula for the 13 primes up to 41 with $\mu(X) = 1.48450$, which appears
to have been chosen subject to the constraint $\max(X) < 2^{32}$.
Gourdon and Sebah~\cite{gourdon2004logarithmic} give a much less efficient
formula for the first 25 primes derived from \eqref{eq:logpdiff}, with $\mu(X) > 7.45186$.

The claim that the formulas in Tables~\ref{tab:logrelations} and ~\ref{tab:atanrelations}
are optimal comes with several caveats.
We can achieve lower Lehmer measures
if we add more arctangents.
Indeed, the formula for $P = \{2, 3, 5, 7\}$ has a
lower Lehmer measure than the formulas for $\{2\}$, $\{2, 3\}$ and $\{2, 3, 5\}$,
so we may just as well compute four logarithms if we want the first one or three.
A more efficient formula
for $\log(2)$ alone is the three-term $X = \{26, 4801, 8749\}$ with $\mu(X) = 1.232 05$
which
however cannot be used to compute $\log(3)$, $\log(5)$ or $\log(7)$
(the set $X^2-1$ is 7-smooth but does not yield a relation for either 3, 5 or 7).
The 1-term formula for $\atan(1) = \pi/4$ has infinite
Lehmer measure
while Machin's formula \eqref{eq:machin},
which follows from the 13-smooth factorizations
$5^2 + 1 = 2 \cdot 13$ and $239^2+1 = 2 \cdot 13^4$,
achieves $\mu(X) = 1.85112$.

In practice $\mu(X)$ is not necessarily an accurate measure of efficiency: it overestimates the benefits
of increasing $x$, essentially because the running
time in binary splitting tends to be dominated by the top-level
multiplications which are independent of the number of leaf nodes.
It is therefore likely an advantage to keep the number
of arctangents close to $n$.


A curiosity is that in the logarithm relations,
we have $\det(R) = \pm 1$ and therefore $M \in \ZZ^{n \times n}$
for the first 21 sets of primes $P$, but for
$P$ containing the primes up to 79, 83, 89 and 97 respectively
the determinants are $-2$, $-6$, $-4$ and $-4$.

\begin{table}
\centering
\caption{\small $n$-term Machin formulas $\{\operatorname{atanh}(1/x) : x \in X\}$ for simultaneous computation of $\log(p)$ for the first $n$ primes $p \in P$.}
\label{tab:logrelations}
\tiny
\renewcommand{\arraystretch}{1.2}
\begin{tabular}{ r | l | >{\raggedright}p{9cm} | l }
$n$ & $P$ & $X$ & $\mu(X)$ \\
\hline
1 & 2 & 3 & 2.09590 \\
2 & 2, 3 & 7, 17 & 1.99601 \\
3 & 2, 3, 5 & 31, 49, 161 & 1.71531 \\
4 & 2 \ldots 7 & 251, 449, 4801, 8749 & 1.31908 \\
5 & 2 \ldots 11 & 351, 1079, 4801, 8749, 19601 & 1.48088 \\
6 & 2 \ldots 13 & 1574, 4801, 8749, 13311, 21295, 246401 & 1.49710 \\
7 & 2 \ldots 17 & 8749, 21295, 24751, 28799, 74359, 388961, 672281 & 1.49235 \\
8 & 2 \ldots 19 & 57799, 74359, 87361, 388961, 672281, 1419263, 11819521, 23718421 & 1.40768 \\
9 & 2 \ldots 23 & 143749, 672281, 1419263, 1447874, 4046849, 8193151, 10285001, 11819521, 23718421 & 1.40594 \\
10 & 2 \ldots 29 & 1419263, 1447874, 11819521, 12901780, 16537599, 23718421, 26578124, 36171409, 192119201, 354365441 & 1.38570 \\
11 & 2 \ldots 31 & 1447874, 11819521, 12901780, 16537599, 23718421, 36171409, 287080366, 354365441, 362074049, 740512499, 3222617399 & 1.42073 \\
12 & 2 \ldots 37 & 36171409, 42772001, 55989361, 100962049, 143687501, 287080366, 362074049, 617831551, 740512499, 3222617399, 6926399999, 9447152318 & 1.40854 \\
13 & 2 \ldots 41 & 51744295, 170918749, 265326335, 287080366, 362074049, 587270881, 831409151, 2470954914, 3222617399, 6926399999, 9447152318, 90211378321, 127855050751 & 1.42585 \\
14 & 2 \ldots 43 & 287080366, 975061723, 980291467, 1181631186, 1317662501, 2470954914, 3222617399, 6926399999, 9447152318, 22429958849, 36368505601, 90211378321, 127855050751, 842277599279 & 1.43055 \\
15 & 2 \ldots 47 & 2470954914, 2473686799, 3222617399, 4768304960, 6926399999, 9447152318, 22429958849, 36974504449, 74120970241, 90211378321, 127855050751, 384918250001, 569165414399, 842277599279, 2218993446251 & 1.42407 \\
16 & 2 \ldots 53 & 9943658495, 15913962107, 19030755899, 22429958849, 22623739319, 36974504449, 90211378321, 123679505951, 127855050751, 187753824257, 384918250001, 569165414399, 842277599279, 1068652740673, 2218993446251, 2907159732049 & 1.44292 \\
17 & 2 \ldots 59 & 22429958849, 56136455649, 92736533231, 122187528126, 123679505951, 127855050751, 134500454243, 187753824257, 384918250001, 569165414399, 842277599279, 1829589379201, 2218993446251, 2569459276099, 2907159732049, 22518692773919, 41257182408961 & 1.45670 \\
18 & 2 \ldots 61 & 123679505951, 210531506249, 367668121249, 384918250001, 711571138431, 842277599279, 1191139875199, 1233008445689, 1829589379201, 2218993446251, 2569459276099, 2907159732049, 3706030044289, 7233275252995, 9164582675249, 22518692773919, 41257182408961, 63774701665793 & 1.46360 \\
19 & 2 \ldots 67 & 664954699135, 842277599279, 932784765626, 1191139875199, 1233008445689, 1726341174999, 1829589379201, 2198699269535, 2218993446251, 2569459276099, 2907159732049, 3706030044289, 7233275252995, 8152552404881, 9164582675249, 22518692773919, 25640240468751, 41257182408961, 63774701665793 & 1.51088 \\
20 & 2 \ldots 71 & 932784765626, 1986251708497, 2200009162625, 2218993446251, 2907159732049, 5175027061249, 7233275252995, 8152552404881, 8949772845287, 9164582675249, 12066279000049, 13055714577751, 22518692773919, 25640240468751, 31041668486401, 41257182408961, 63774701665793, 115445619421397, 121336489966251, 238178082107393 & 1.52917 \\
21 & 2 \ldots 73 & 7233275252995, 8152552404881, 8949772845287, 9164582675249, 10644673332721, 13055714577751, 21691443063179, 22518692773919, 25640240468751, 25729909301249, 41257182408961, 54372220771987, 63774701665793, 103901723427151, 106078311729181, 114060765404951, 115445619421397, 121336489966251, 238178082107393, 1796745215731101, 4573663454608289 & 1.53515 \\
22 & 2 \ldots 79 & 38879778893521, 41257182408961, 44299089391103, 62678512919879, 63774701665793, 69319674756179, 70937717129551, 103901723427151, 106078311729181, 114060765404951, 115445619421397, 117774370786951, 121336489966251, 217172824950401, 238178082107393, 259476225058051, 386624124661501, 478877529936961, 1796745215731101, 2767427997467797, 4573663454608289, 19182937474703818751 & 1.52802 \\
23 & 2 \ldots 83 & 103901723427151, 112877019076249, 114060765404951, 115445619421397, 117774370786951, 121336489966251, 134543112911873, 148569359956291, 201842423186689, 206315395261249, 217172824950401, 238178082107393, 259476225058051, 386624124661501, 473599589105798, 478877529936961, 1796745215731101, 1814660314218751, 2767427997467797, 4573663454608289, 17431549081705001, 34903240221563713, 19182937474703818751 & 1.55501 \\
24 & 2 \ldots 89 & 134543112911873, 148569359956291, 166019820559361, 201842423186689, 206315395261249, 211089142289024, 217172824950401, 238178082107393, 259476225058051, 330190746672799, 386624124661501, 473599589105798, 478877529936961, 1796745215731101, 1814660314218751, 2767427997467797, 2838712971108351, 4573663454608289, 9747977591754401, 11305332448031249, 17431549081705001, 34903240221563713, 332110803172167361, 19182937474703818751 & 1.58381 \\
25 & 2 \ldots 97 & 373632043520429, 386624124661501, 473599589105798, 478877529936961, 523367485875499, 543267330048757, 666173153712219, 1433006524150291, 1447605165402271, 1744315135589377, 1796745215731101, 1814660314218751, 2236100361188849, 2767427997467797, 2838712971108351, 3729784979457601, 4573663454608289, 9747977591754401, 11305332448031249, 17431549081705001, 21866103101518721, 34903240221563713, 99913980938200001, 332110803172167361, 19182937474703818751 & 1.60385 \\
\end{tabular}
\end{table}


\begin{table}
\centering
\caption{\small $n$-term Machin formulas $\{\operatorname{atan}(1/x) : x \in X\}$ for simultaneous computation of the irreducible angles
$\operatorname{atan}(b/a)$ for the first $n$ nonreal Gaussian primes $a+bi$, having norms $a^2+b^2 \in Q$.}
\label{tab:atanrelations}
\tiny
\renewcommand{\arraystretch}{1.2}
\begin{tabular}{ r | l | >{\raggedright}p{9cm} | l }
$n$ & $Q$ & $X$ & $\mu(X)$ \\
\hline
1 & 2 & 1 & $\infty$ \\
2 & 2, 5 & 3, 7 & 3.27920 \\
3 & 2, 5, 13 & 18, 57, 239 & 1.78661 \\
4 & 2 \ldots 17 & 38, 57, 239, 268 & 2.03480 \\
5 & 2 \ldots 29 & 38, 157, 239, 268, 307 & 2.32275 \\
6 & 2 \ldots 37 & 239, 268, 307, 327, 882, 18543 & 2.20584 \\
7 & 2 \ldots 41 & 268, 378, 829, 882, 993, 2943, 18543 & 2.33820 \\
8 & 2 \ldots 53 & 931, 1772, 2943, 6118, 34208, 44179, 85353, 485298 & 2.01152 \\
9 & 2 \ldots 61 & 5257, 9466, 12943, 34208, 44179, 85353, 114669, 330182, 485298 & 1.95679 \\
10 & 2 \ldots 73 & 9466, 34208, 44179, 48737, 72662, 85353, 114669, 330182, 478707, 485298 & 2.03991 \\
11 & 2 \ldots 89 & 51387, 72662, 85353, 99557, 114669, 157318, 260359, 330182, 478707, 485298, 24208144 & 2.06413 \\
12 & 2 \ldots 97 & 157318, 330182, 390112, 478707, 485298, 617427, 1984933, 2343692, 3449051, 6225244, 22709274, 24208144 & 1.96439 \\
13 & 2 \ldots 101 & 683982, 1984933, 2343692, 2809305, 3014557, 6225244, 6367252, 18975991, 22709274, 24208144, 193788912, 201229582, 2189376182 & 1.84765 \\
14 & 2 \ldots 109 & 2298668, 2343692, 2809305, 3014557, 6225244, 6367252, 18975991, 22709274, 24208144, 168623905, 193788912, 201229582, 284862638, 2189376182 & 1.91451 \\
15 & 2 \ldots 113 & 2343692, 2809305, 3801448, 6225244, 6367252, 7691443, 18975991, 22709274, 24208144, 168623905, 193788912, 201229582, 284862638, 599832943, 2189376182 & 2.01409 \\
16 & 2 \ldots 137 & 4079486, 6367252, 7691443, 8296072, 9639557, 10292025, 18975991, 19696179, 22709274, 24208144, 168623905, 193788912, 201229582, 284862638, 599832943, 2189376182 & 2.12155 \\
17 & 2 \ldots 149 & 9689961, 10292025, 13850847, 18975991, 19696179, 22709274, 24208144, 32944452, 58305593, 60033932, 168623905, 193788912, 201229582, 284862638, 314198789, 599832943, 2189376182 & 2.18157 \\
18 & 2 \ldots 157 & 22709274, 32944452, 58305593, 60033932, 127832882, 160007778, 168623905, 193788912, 201229582, 284862638, 299252491, 314198789, 361632045, 599832943, 851387893, 2189376182, 2701984943, 3558066693 & 2.14866 \\
19 & 2 \ldots 173 & 127832882, 160007778, 168623905, 193788912, 201229582, 299252491, 314198789, 327012132, 361632045, 599832943, 851387893, 1117839407, 2189376182, 2701984943, 3558066693, 12139595709, 12957904393, 120563046313, 69971515635443 & 2.09258 \\
20 & 2 \ldots 181 & 299252491, 314198789, 327012132, 361632045, 599832943, 851387893, 1112115023, 1117839407, 1892369318, 2189376182, 2701984943, 2971354082, 3558066693, 5271470807, 12139595709, 12957904393, 14033378718, 18986886768, 120563046313, 69971515635443 & 2.10729 \\
21 & 2 \ldots 193 & 1112115023, 1117839407, 1479406293, 1696770582, 1892369318, 2112819717, 2189376182, 2701984943, 2971354082, 3558066693, 4038832337, 5271470807, 7959681215, 8193535810, 12139595709, 12957904393, 14033378718, 18710140581, 18986886768, 120563046313, 69971515635443 & 2.13939 \\
22 & 2 \ldots 197 & 1479406293, 1892369318, 2112819717, 2189376182, 2701984943, 2971354082, 3558066693, 4038832337, 5271470807, 6829998457, 7959681215, 8193535810, 12139595709, 12185104420, 12957904393, 14033378718, 18710140581, 18986886768, 20746901917, 104279454193, 120563046313, 69971515635443 & 2.19850 \\
\end{tabular}
\end{table}

\section{Implementation results}

\label{sect:implementation}

The algorithms have been implemented in Arb~\cite{Joh2017} version 2.23.
The following results
were obtained with Arb~2.23 linked against GMP~6.2.1~\cite{GMP}, MPFR~4.1.0~\cite{Fousse2007}, and FLINT~2.9~\cite{Hart2010},
running on an AMD Ryzen 7 PRO 5850U (Zen3).

\subsection{Default implementations with fixed $n$}

Previously, all elementary functions in Arb used Taylor series
with precomputed lookup tables up to $B = 4608$ bits.
The tables are $m$-partite giving $r$-bit reduction
with $r \le 14$ and $m \le 2$,
requiring 236 KB of fixed storage~\cite{Johansson2015elementary}.
At higher precision, the previous implementations used argument reduction based on
repeated argument-halving (requiring squaring or square roots) together with rectangular splitting or
bit-burst evaluation of Taylor series,
with the exception of log which wrapped the AGM-based logarithm in MPFR.
To the author's knowledge, these were the fastest arbitrary-precision
implementations of elementary functions available in
public software libraries prior to this work.

In Arb 2.23, all the elementary functions were rewritten
to use the new algorithm
with the fixed number $n = 13$ of primes,
starting from a precision between $B = 2240$ bits (for exp) and
$B = 3400$ bits (for atan) up to $B = 4000000$ bits (just over one million digits).
The Newton iterations~\eqref{eq:logy} and \eqref{eq:atany} are used to reduce log and atan to the
exponential and trigonometric functions.
The $B$-bit precomputations
of logarithms and arctangents are done at runtime using the $n = 13$
Machin-like formulas of Table~\ref{tab:logrelations} and Table~\ref{tab:atanrelations}.

We compare timings for the old and new implementations in~Table \ref{tab:oldnew13}.

\begin{table}
\centering
\caption{\small Time to compute elementary functions to $D$ decimal digits ($B \approx 3.32 D$) with Arb 2.23. \emph{Old} is the time in seconds with the new algorithm disabled. \emph{New} is the time in seconds
with the new algorithm enabled, using the fixed default number $n = 13$ of primes. \emph{First} is the time for a first function call, and \emph{Repeat} is the time for repeated calls (with logarithms and
other data already cached). We show average timings for 100 uniformly random input $x \in (0, 2)$.}
\label{tab:oldnew13}
\scriptsize
\setlength{\tabcolsep}{3pt}
\renewcommand{\arraystretch}{1.4}
\begin{tabular}{ c c | c c | c c | c c | c c }
\multicolumn{2}{c}{}  & \multicolumn{2}{c}{$\exp(x)$} & \multicolumn{2}{c}{$\log(x)$} & \multicolumn{2}{c}{$(\cos(x), \sin(x))$} & \multicolumn{2}{c}{$\atan(x)$} \\ \hline
$D$ & & First & Repeat & First & Repeat & First & Repeat & First & Repeat \\ \hline
1000    & Old & 2.92e-05 & 2.91e-05 & 0.000145 & 3.69e-05 & 3.49e-05 & 3.49e-05 & 3.52e-05 & 3.52e-05 \\
        & New & 0.000182 & 2.04e-05 & 0.000188 & 2.58e-05 & 0.00019 & 2.84e-05 & 3.52e-05 & 3.52e-05 \\
        & Speedup & 0.16$\times$ & 1.43$\times$ & 0.77$\times$ & 1.43$\times$ & 0.18$\times$ & 1.23$\times$ & 1.00$\times$ & 1.00$\times$ \\ \hline
2000    & Old & 0.000103 & 0.000101 & 0.000367 & 0.000110 & 0.000217 & 9.92e-05 & 0.000423 & 0.000217 \\
        & New & 0.000480 & 4.9e-05 & 0.000500 & 6.07e-05 & 0.000542 & 7.92e-05 & 0.000564 & 9.83e-05 \\
        & Speedup & 0.22$\times$ & 2.06$\times$ & 0.73$\times$ & 1.81$\times$ & 0.40$\times$ & 1.25$\times$ & 0.75$\times$ & 2.21$\times$ \\ \hline
4000    & Old & 0.000355 & 0.000353 & 0.00103 & 0.000348 & 0.000511 & 0.000341 & 0.000915 & 0.000660 \\
        & New & 0.00107 & 0.000149 & 0.00111 & 0.000187 & 0.00119 & 0.000211 & 0.00124 & 0.000269 \\
        & Speedup & 0.33$\times$ & 2.37$\times$ & 0.93$\times$ & 1.86$\times$ & 0.43$\times$ & 1.62$\times$ & 0.74$\times$ & 2.45$\times$ \\ \hline
10000   & Old & 0.00185 & 0.00168 & 0.00439 & 0.00166 & 0.0022 & 0.00177 & 0.00323 & 0.00272 \\
        & New & 0.00384 & 0.000826 & 0.00418 & 0.000977 & 0.00417 & 0.000935 & 0.00461 & 0.00122 \\
        & Speedup & 0.48$\times$ & 2.03$\times$ & 1.05$\times$ & 1.70$\times$ & 0.53$\times$ & 1.89$\times$ & 0.70$\times$ & 2.23$\times$ \\ \hline
100000  & Old & 0.0541 & 0.0536 & 0.143 & 0.0632 & 0.0880 & 0.0818 & 0.0957 & 0.0896 \\
        & New & 0.107 & 0.0354 & 0.114 & 0.0377 & 0.129 & 0.0509 & 0.140 & 0.0586 \\
        & Speedup & 0.51$\times$ & 1.52$\times$ & 1.25$\times$ & 1.68$\times$ & 0.68$\times$ & 1.61$\times$ & 0.68$\times$ & 1.53$\times$ \\ \hline
1000000 & Old & 1.10 & 1.09 & 2.84 & 1.36 & 1.66 & 1.61 & 2.02 & 1.97 \\
        & New & 2.18 & 0.864 & 2.31 & 0.982 & 2.83 & 1.25 & 3.02 & 1.58 \\
        & Speedup & 0.51$\times$ & 1.26$\times$ & 1.23$\times$ & 1.39$\times$ & 0.59$\times$ & 1.29$\times$ & 0.67$\times$ & 1.25$\times$ \\
\end{tabular}
\end{table}

\subsubsection*{Remarks}

The average speedup is around a factor two ($1.3\times$ to $2.4\times$) over a large range of precisions.
The typical slowdown for a first function call is also roughly a factor two, i.e.\
the precomputation takes about as long as a single extra function call.\footnote{The figures are a bit worse at lower precision due to various overheads which could be avoided.}
This
is clearly a worthwhile tradeoff for most applications;
e.g.\ for a numerical integration $\smash \int_a^b f(x) dx$
where the integrand $f$ will be evaluated many times,
we do observe a factor-two speedup
in the relevant precision ranges.

The relatively large speedup for atan is explained by the fact that the traditional argument reduction method involves repeated square roots which are a significant constant factor more expensive than the squarings for exp.

The relatively small speedup for sin and cos is explained by the fact that traditional argument reduction method only requires real squarings (via the half-angle formula for cos), while the new method uses complex arithmetic.

Previously, the AGM-based logarithm was neck and neck with the Taylor series for exp at any precision (these algorithms were therefore roughly interchangeable if one were to use Newton iteration to compute one function from the other). With the new algorithm, Taylor series have a clear lead.

The default parameter $n = 13$ was chosen to optimize performance
around a few thousand digits, this range being more important for typical
applications than millions of digits.
As shown below, it is possible to achieve larger speedup at very high
precision by choosing a larger $n$.

\subsection{Precomputation of reduction tables}

Table~\ref{tab:statictime} shows sample results
for the precomputation phase of Algorithm~\ref{alg:linred}
to generate tables of approximate relations over $n$ logarithms
or arctangents.

\begin{table}
\centering
\caption{\small Static precomputation of reduction tables: phase (1) of Algorithm~\ref{alg:linred}.}
\label{tab:statictime}
\scriptsize
\renewcommand{\arraystretch}{1.4}
\begin{tabular}{ c c c c c c }
$\alpha_1, \ldots, \alpha_n$ & $n$  & Smallest $\varepsilon_i$ & Max $r$ & Data & Time \\ \hline
Logarithms & 2    & $\varepsilon_{7} = +1.82 \cdot 10^{-5}$ & 15 & 0.2 KiB & 0.0000514 s \\
& 4    & $\varepsilon_{11} = -1.46 \cdot 10^{-14}$ & 45 & 0.3 KiB & 0.000228 s \\
& 8    & $\varepsilon_{33} = +7.66 \cdot 10^{-33}$ & 106 & 1.1 KiB & 0.00249 s \\
& 16   & $\varepsilon_{67} = +5.18 \cdot 10^{-71}$ & 233 & 3.2 KiB & 0.0447 s \\
& 32   & $\varepsilon_{144} = -1.51 \cdot 10^{-141}$ & 467 & 11 KiB & 1.24 s \\
& 64   & $\varepsilon_{268} = -4.42 \cdot 10^{-266}$ & 881 & 38 KiB & 34.2 s \\ \hline
Arctangents & 2   & $\varepsilon_{7} = -4.75 \cdot 10^{-5}$ & 14 & 0.2 KiB & 0.0000472 s \\
& 4   & $\varepsilon_{14} = -2.95 \cdot 10^{-15}$ & 48 & 0.4 KiB & 0.000248 s \\
& 8   & $\varepsilon_{33} = +6.43 \cdot 10^{-33}$ & 106 & 1.1 KiB & 0.00256 s \\
& 16   & $\varepsilon_{64} = +1.77 \cdot 10^{-71}$ & 235 & 3.0 KiB & 0.0448 s \\
& 32   & $\varepsilon_{143} = +1.70 \cdot 10^{-140}$ & 464 & 11 KiB & 1.22 s \\
& 64   & $\varepsilon_{270} = +1.42 \cdot 10^{-267}$ & 886 & 38 KiB & 34.6 s \\
\end{tabular}
\end{table}

Here we choose the convergence factor $C = 10$ (each approximate relation $\varepsilon_i$ adds one decimal)
and we terminate before the first relation with
a coefficient $|d_{i,j}| \ge 2^{15}$. This bound was chosen
for convenience of storing
table entries in 16-bit integers;
it is also a reasonable cutoff since
larger exponents will pay off only for multi-million $B$ (as we will see below).
We test the method up to $n = 64$, where the smallest tabulated $\varepsilon_i$ corresponds to
an argument reduction of more than $r = 800$ bits.\footnote{Part of the implementation uses machine-precision
floating-point numbers with a limited exponent range,
making $|\varepsilon_i| < 2^{-1024} \approx 10^{-300}$ inaccessible.
Like the 16-bit limit, this is again a trivial technical restriction
which we do not bother to lift
since there would be a pay-off only for multi-million $B$.}

Since the tables are small (a few KiB) and independent of $B$, they
can be precomputed once and for all,
so the timings (here essentially just exercising FLINT's LLL implementation)
are not really relevant.
Indeed, in the previously
discussed default implementation of elementary functions,
the $n = 13$ tables are stored as static arrays
written down in the source code.
However, the timings are reasonable enough that
tables could be generated at runtime in applications
that will perform a large number of function evaluations.

\begin{table}
\centering
\caption{\small Computation of the exponential function and the trigonometric functions. The argument is taken to be $x = \sqrt{2}-1$.
\emph{Precomp} is the time (in seconds) to precompute $n$ logarithms or arctangents for use at $B$-bit precision. The cached logarithms or arctangents take up \emph{Data} space.
\emph{Time} is the time to evaluate the function once this data has been precomputed. The argument is reduced to size $2^{-r}$.}
\label{tab:functime}
\scriptsize
\renewcommand{\arraystretch}{1.2}
\begin{tabular}{ c c c | l c l | l c l }
\multicolumn{3}{c}{}  & \multicolumn{3}{c}{$\exp(x)$} & \multicolumn{3}{c}{$\cos(x) + i \sin(x) = \exp(ix)$} \\
$B$ & $n$ & Data & Precomp & $r$ & Time & Precomp & $r$ & Time \\ \hline
3333  &  0   &  &     &    &  2.89e-05  &   &   & 3.56e-05 \\
  & 2  &  0.8 KiB  &  5.33e-05  &  11  &  2.88e-05  &         6.34e-05  &  11  &  3.49e-05 \\
  & 4  &  1.6 KiB  &  5.42e-05  &  15  &  2.71e-05  &         7.35e-05  &  22  &  2.74e-05 \\
  & 8  &  3.3 KiB  &  7.61e-05  &  32  &  2.06e-05  &         9.65e-05  &  33  &  2.72e-05 \\
  & 16  &  6.5 KiB  &  0.000131  &  73  &  1.78e-05  &         0.000136  &  37  &  2.89e-05 \\
  & 32  &  13.0 KiB  &  0.000268  &  60  &  1.97e-05  &         0.000411  &  38  &  2.92e-05 \\
  & 64  &  26.0 KiB  &  0.000605  &  60  &  2.2e-05  &         0.00104  &  38  &  3.15e-05 \\ \hline
10000  &  0   &  &     &    &  0.000202  &   &   & 0.000207 \\
  & 2  &  2.4 KiB  &  0.000238  &  11  &  0.000183  &         0.000281  &  13  &  0.000209 \\
  & 4  &  4.9 KiB  &  0.000240  &  27  &  0.000137  &         0.000333  &  30  &  0.000159 \\
  & 8  &  9.8 KiB  &  0.000335  &  52  &  0.000106  &         0.000412  &  41  &  0.000144 \\
  & 16  &  19.5 KiB  &  0.000579  &  83  &  8.48e-05  &         0.000633  &  61  &  0.000114 \\
  & 32  &  39.1 KiB  &  0.00123  &  86  &  8.75e-05  &         0.00187  &  47  &  0.000129 \\
  & 64  &  78.1 KiB  &  0.00270  &  72  &  9.71e-05  &         0.00468  &  47  &  0.000131 \\ \hline
33333  &  0   &  &     &    &  0.00166  &   &   & 0.00178 \\
  & 2  &  8.1 KiB  &  0.00135  &  18  &  0.00135  &         0.0016  &  13  &  0.00167 \\
  & 4  &  16.3 KiB  &  0.00136  &  44  &  0.00107  &         0.00186  &  30  &  0.00133 \\
  & 8  &  32.6 KiB  &  0.00199  &  56  &  0.000938  &         0.00239  &  65  &  0.00110 \\
  & 16  &  65.1 KiB  &  0.00330  &  89  &  0.000748  &         0.00371  &  90  &  0.000932 \\
  & 32  &  130.2 KiB  &  0.00683  &  139  &  0.000637  &         0.0103  &  138  &  0.000841 \\
  & 64  &  260.4 KiB  &  0.0152  &  168  &  0.000614  &         0.0256  &  63  &  0.00103 \\ \hline
100000  &  0   &  &     &    &  0.00895  &   &   & 0.0125 \\
  & 2  &  24.4 KiB  &  0.00679  &  18  &  0.00747  &         0.00786  &  17  &  0.0119 \\
  & 4  &  48.8 KiB  &  0.0068  &  44  &  0.00638  &         0.00922  &  40  &  0.00987 \\
  & 8  &  97.7 KiB  &  0.00977  &  71  &  0.00565  &         0.0119  &  65  &  0.00754 \\
  & 16  &  195.3 KiB  &  0.0164  &  106  &  0.00534  &         0.0179  &  90  &  0.00625 \\
  & 32  &  390.6 KiB  &  0.0337  &  161  &  0.00445  &         0.0491  &  138  &  0.00523 \\
  & 64  &  781.2 KiB  &  0.0755  &  240  &  0.00383  &         0.125  &  126  &  0.00612 \\ \hline
1000000  &  0   &  &     &    &  0.221  &   &   & 0.337 \\
  & 2  &  244.1 KiB  &  0.159  &  18  &  0.195  &         0.187  &  17  &  0.322 \\
  & 4  &  488.3 KiB  &  0.159  &  47  &  0.175  &         0.219  &  40  &  0.295 \\
  & 8  &  976.6 KiB  &  0.228  &  99  &  0.154  &         0.271  &  96  &  0.273 \\
  & 16  &  1.9 MiB  &  0.37  &  142  &  0.140  &         0.419  &  118  &  0.260 \\
  & 32  &  3.8 MiB  &  0.77  &  161  &  0.136  &         1.14  &  171  &  0.255 \\
  & 64  &  7.6 MiB  &  1.72  &  454  &  0.120  &         2.91  &  391  &  0.178 \\ \hline
10000000  &  0   &  &     &    &  4.36  &   &   & 6.50 \\
  & 2  &  2.4 MiB  &  3.02  &  18  &  3.89  &         3.56  &  17  &  6.18 \\
  & 4  &  4.8 MiB  &  3.01  &  47  &  3.53  &         4.1  &  40  &  5.75 \\
  & 8  &  9.5 MiB  &  4.14  &  110  &  3.18  &         5.03  &  109  &  5.24 \\
  & 16  &  19.1 MiB  &  6.57  &  222  &  2.90  &         7.49  &  203  &  5.13 \\
  & 32  &  38.1 MiB  &  13.8  &  338  &  2.61  &         20.6  &  348  &  4.64 \\
  & 64  &  76.3 MiB  &  31.3  &  551  &  2.39  &         53.4  &  592  &  4.50 \\
\end{tabular}
\end{table}

\subsection{Function evaluation with variable $n$ and $B$}

Table~\ref{tab:functime} shows timings for the computation
of the exponential function and trigonometric functions
for different combinations of precision $B$ and number of primes $n$.
The $n = 0$ reference timings correspond
to the old algorithm without precomputation,
in which repeated squaring will be used instead.

At lower precisions, using 10-20 primes seems to be optimal.
It is interesting to note that roughly a factor-two speedup
can be achieved across a huge range of precision
when $n$ increases with $B$.
It seems likely that $n = 128$ or more primes could
be useful at extreme precision,
though the precomputation will increase proportionally.

We used the Machin-like formulas from Table~\ref{tab:logrelations} and Table~\ref{tab:atanrelations}
only up to $n = 25$ or $n = 22$;
for $n = 32$ and $n = 64$ we fall back
on less optimized formulas,
which results in a noticeably slower precomputation.

\section{Extensions and further work}

We conclude with some ideas for future research.

\subsection{Complexity analysis and fine-tuning}

It would be interesting to perform a more precise
complexity analysis.
Under some assumptions about the underlying
arithmetic, it should be possible to obtain a theoretical prediction
for the optimal number of primes $n$ as a function of the bit precision $B$,
with an estimate of the possible speedup when $B \to \infty$.

There are a large number of free parameters in the
new algorithm (the number of primes $n$,
the choice of primes,
the precise setup of
the precomputed relations $\varepsilon_i$,
the allowed size of the power product,
choices in the subsequent Taylor series evaluation...).
Timings can fluctuate depending on
small adjustments to these parameters
and with different values of the argument~$x$.
It is plausible that a consistent speedup
can be obtained by simply tuning all the parameters more
carefully.

\subsection{Complex arguments}

All elementary functions of complex arguments can be decomposed
into real exponentials, logarithms, sines and arctangents
after separating real and imaginary parts.
An interesting alternative would be to
compute $\exp(z)$ or $\log(z)$
directly over $\mathbb{C}$,
reducing $z$ with respect to complex lattices
generated by pairs of Gaussian primes.
We do not know whether this presents any advantages
over separating the components.

\subsection{$p$-adic numbers}

The same methods should work in the $p$-adic setting.
For the $p$-adic exponential and logarithm,
we can choose a basis of $n$ prime numbers $p_i \ne p$
and use LLL to precompute relations
$\sum_{i=1}^n c_i \log(p_i) = O(p^i)$
for $i = 1, 2, \ldots, r$.
We can then use these relations to reduce
the argument to order $O(p^r)$
before evaluating the function
using Taylor series or the $p$-adic bit-burst method~\cite{caruso2021fast}.
We have not attempted to analyze
or implement this algorithm.

\subsection{More Machin-like formulas}

It would be useful to have larger tables of optimized Machin-like
formulas for multi-evaluation of logarithms and arctangents.
In practice, formulas need not be optimal as long as they
are ``good enough''; for example, one could restrict the search
space to 64-bit arctangent arguments $x$.
Nevertheless, a large-scale computation
of theoretically optimal tables would be an interesting
challenge of its own.

\section{Acknowledgements}

The author was present at RISC in 2011 where
Arnold Sch\"{o}nhage gave one of the talks~\cite{schonhage2011} presenting
his original ``medium-precision'' version of the algorithm
using a pair of primes.
Ironically, the author has no memory of the event
beyond the published talk abstract;
the inspiration for the present work came much later,
with Machin-like formulas for logarithms as the starting point,
and the details herein were developed independently.
Nevertheless, Sch\"{o}nhage certainly deserves credit
for the core idea.
We have tried unsuccessfully to contact Sch\"{o}nhage (who is now retired)
for notes about his version
of the algorithm.

The author learned about the process to
find Machin-like formulas
thanks to MathOverflow comments
by Douglas Zare and the user ``Meij''~\cite{125687}
explaining the method and pointing to the relevant chapter in Arndt's book.

Algorithm~\ref{alg:linred} was inspired by a comment by
Simon Puchert in 2018 proposing an iterative
argument reduction using smooth fractions of the form $(m+1)/m$.
We have substantially improved this algorithm by using
LLL to look for arbitrary smooth fractions
close to 1 instead of restricting to a set of fractions
of special form, and by working with the logarithmic forms
during reduction.

The author was supported by ANR grant ANR-20-CE48-0014-02 NuSCAP.

\bibliographystyle{alpha}
\bibliography{../references}

\end{document}